\documentclass[11pt]{article}

\usepackage{apacite,amsfonts,amsmath,latexsym,graphics,epsfig,subfigure,color,setspace,lineno}
\usepackage{amsfonts,amsmath,amssymb,amsthm,mathtools}
\usepackage[left=3.0cm,right=3.0cm]{geometry}
\usepackage{enumerate,natbib}

\usepackage{bm}

\definecolor{dred}{rgb}{0.5,0,0}
\definecolor{dgreen}{rgb}{0,0.7,0}
\definecolor{dblue}{rgb}{0,0,0.85}
\definecolor{dmagenta}{cmyk}{0,1,0,0.6}
\definecolor{dcyan}{cmyk}{1,0,0,0.5}
\definecolor{grey}{gray}{0.9}
\definecolor{orange}{rgb}{1,0.65,0}
\definecolor{mellow}{rgb}{.847,.72,.525}
\definecolor{golden}{rgb}{.80392,.60784,.11373}
\definecolor{dgolden}{rgb}{.5451,.39608,.03137}
\definecolor{brown}{rgb}{.15,.15,.15}
\definecolor{darkolivegreen}{rgb}{.33333,.41961,.18431}

\def\tr{\mathop{\rm tr}\nolimits}

\def\mR{\mathbb{R}}

\def\cH{{\mathcal H}}

\def\cT{{\mathcal T}}

\def\cX{{\mathcal X}}

\DeclarePairedDelimiter{\norm}{\lVert}{\rVert}

\def\argmin{\mathop{\rm arg\,inf}\limits}

\author{Wicher Bergsma\\London School of Economics and Political Science}
\title{Discussion contribution ``Functional models for time-varying random objects'' by Dubey and M\"uller (to appear in JRSS-B)}

\begin{document}

\maketitle

In an inspiring paper Dubey and M\"uller (DM) extend PCA to the case that observations are metric-valued functions. 
As an alternative, we develop a kernel PCA \citep*[kPCA;][]{ssm98} approach, which we show is closely related to the DM approach. 
While kernel principal components (kPCs) are simply defined, DM require added complexity in the form of ``object FPCs'' and ``Fr\'echet scores''. 

\subsection*{Kernel PCA}

Suppose observations $X_1,\ldots,X_n$ take values in an arbitrary set $\cX$. Let $\cH$ be a Hilbert space and consider an embedding function $e:\cX\rightarrow\cH$. kPCA is essentially PCA on the embedded observations $e(X_1),\ldots,e(X_n)$, whose sample covariance kernel is
\begin{align} \label{cov2}
   C = \frac1n\sum_{i=1}^n (e(X_i)-\overline e)\otimes(e(X_i)-\overline e)\in\cH\otimes\cH.
\end{align}
where $\overline e=n^{-1}\sum e(X_i)$. 
With $\varphi_k$ the $k$th eigenvector of $C$, the $k$th kPC of $x\in\cX$ is 
\begin{align}\label{kpc}
    \big\langle\varphi_k,e(x)\big\rangle_\cH.
\end{align}
The embedding function need not be computed explicitly as the kPCs of the observations are given by the eigenvectors of the $n\times n$ kernel matrix with elements 
\[  K_\cX(X_i,X_j) = \big\langle e(X_i)-\overline e,e(X_j)-\overline e\big\rangle_\cH. \]

\subsection*{Kernel PCA for observations that are metric-valued functions}

If $(\cX,d_\cX)$ is a metric space of negative type, an isometry $e:\cX\rightarrow\cH$ exists, and a corresponding unique centered kernel
\begin{align*}
   K_\cX(x,x') 
   &= \big\langle e(x)-\overline e,e(x')-\overline e\big\rangle_\cH \\
   &= \frac1{2n^2}\sum_{i=1}^{n}\sum_{j=1}^{n}\Big(\norm[\big]{e(x)-e(X_i)}_\cH^2+\norm[\big]{e(x')-e(X_j)}_\cH-\norm[\big]{e(x)-e(x')}_\cH^2-\norm[\big]{e(X_i)^2-e(X_j)}_\cH^2\Big) \\
   &= \frac1{2n^2}\sum_{i=1}^{n}\sum_{j=1}^{n}\Big( d_\cX^2(x,X_i) + d_\cX^2(x',X_j) - d_\cX^2(x,x') - d_\cX^2(X_i,X_j) \Big) .
\end{align*}

Let $(\Omega,d_\Omega)$ be a metric space of negative type and let $\cX$ be the space of $\Omega$-valued functions over an index set $\cT$ for which 
\begin{align*}
   d_\cX^2(x,x') = \int_\cT d_\Omega^2\big(x(t),x'(t)\big)dt.  
\end{align*}
is finite. Then $(\cX,d_\cX)$ is also of negative type, with corresponding kernel
\begin{align*} 
   K_\cX(x,x') 
   &= \frac1{2n^2}\sum_{i=1}^{n}\sum_{j=1}^{n}\Big( d_\cX^2(x,X_i) + d_\cX^2(x',X_j) - d_\cX^2(x,x') - d_\cX^2(X_i,X_j) \Big) \\
   &= \frac12\int_\cT\Big(d_\Omega^2(x(t),X_i(t))+d_\Omega^2(x'(t),X_j(t))-d_\Omega^2(x(t),x'(t))-d_\Omega^2(X_i(t),X_j(t))\Big)dt \nonumber\\
   &= \int_\cT K_\Omega(x(t),x'(t))dt. \label{kx}
\end{align*}

\begin{figure}[tbp]
	\centering
	\includegraphics[width=130mm]{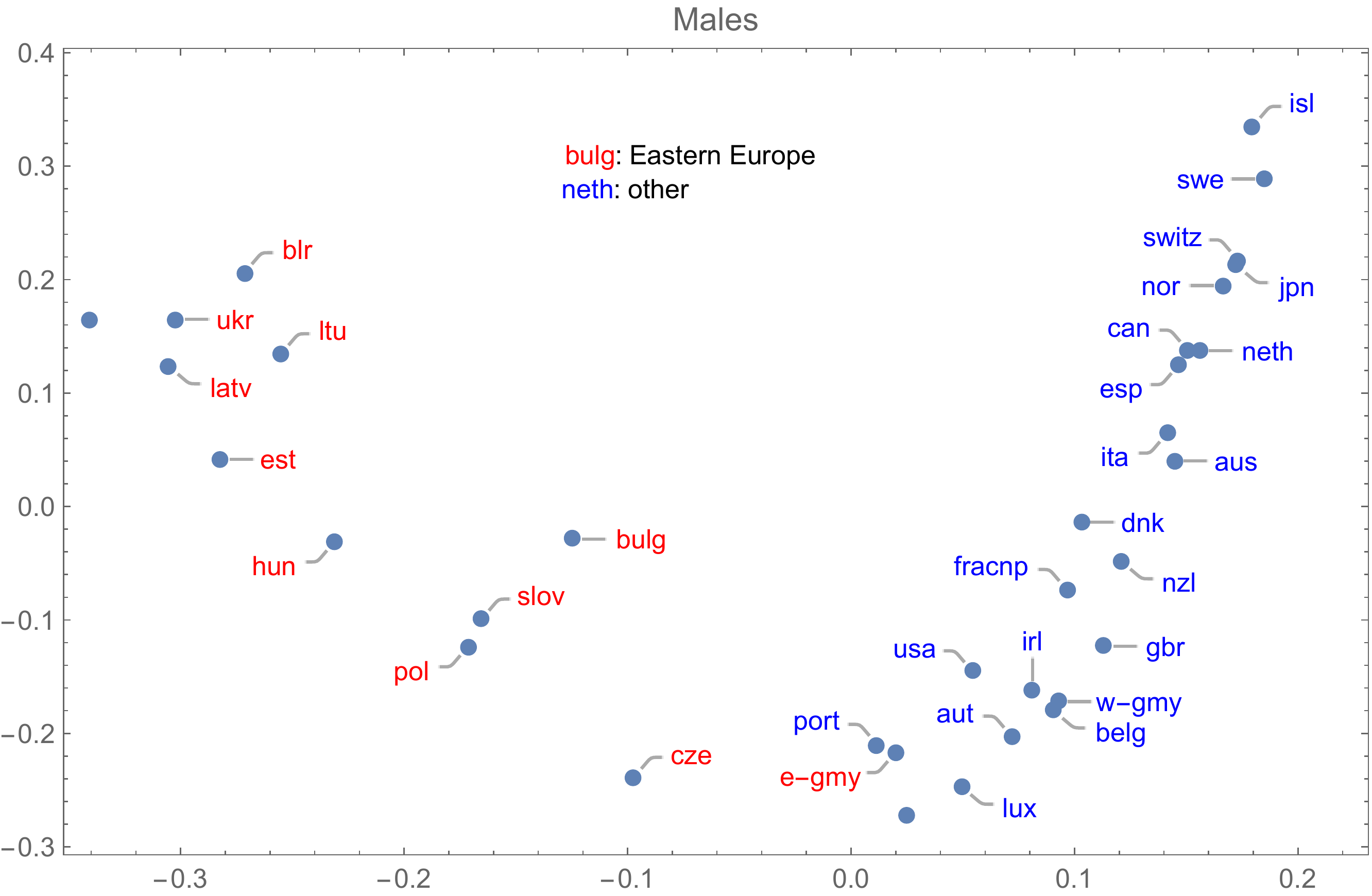} 
	\includegraphics[width=130mm]{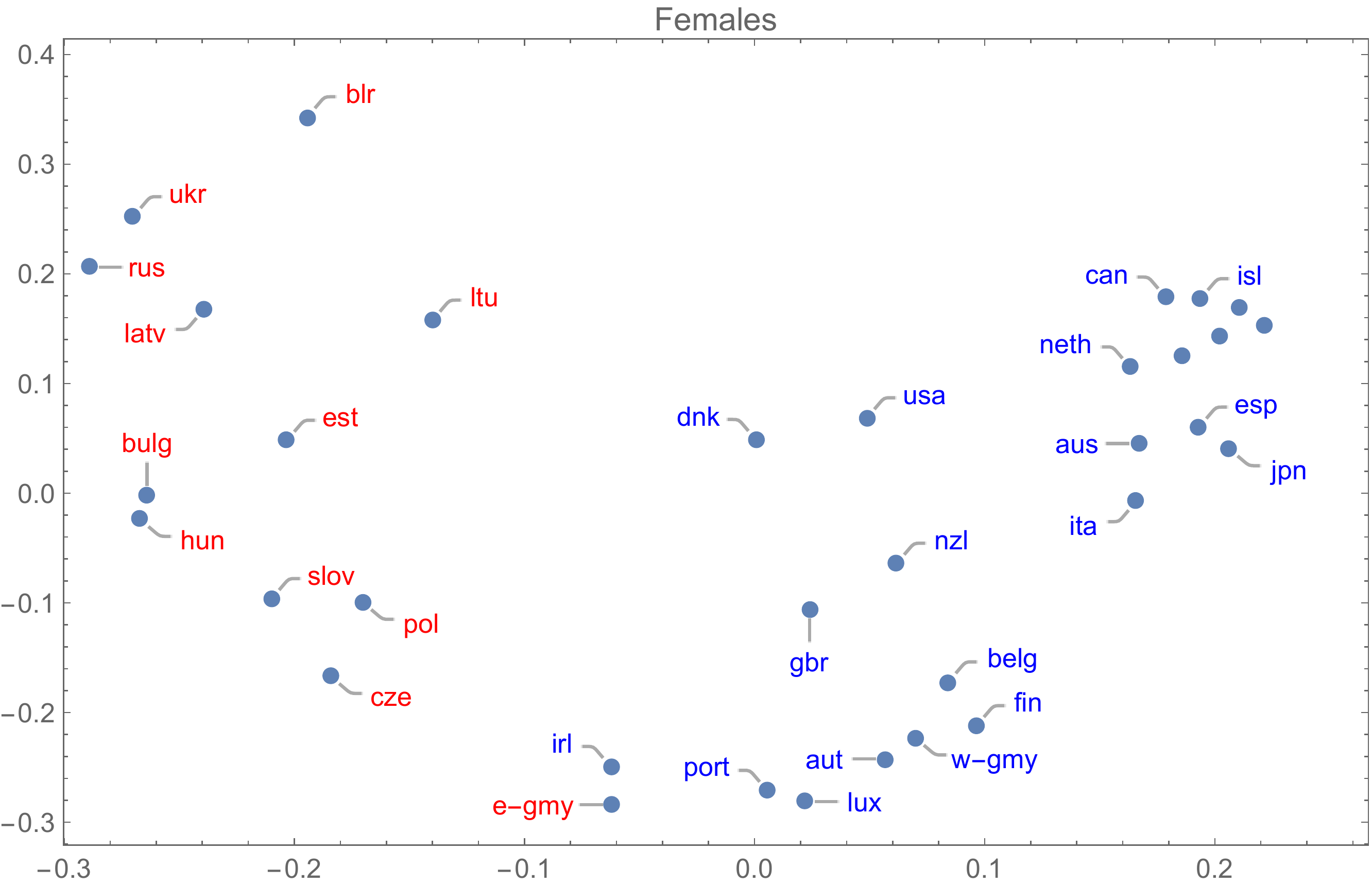}
	\caption{First two kernel principal components for mortality data.}
	\label{fig-pc}
\end{figure}

We reanalysed the mortality data \citep{mortality19} using $L_2$-distance between mortality CDFs. 
Figure~\ref{fig-pc} shows the first two kPCs for both males and females (explained variances around 40\% resp.\ 12\% in both cases). The separation between Eastern European countries and others is clearer than in Figure~10 of DM, which displays Fr\'echet scores. 

\subsection*{Comparison kernel PCA and Dubey-M\"uller method}


Let $e_\Omega:\Omega\rightarrow\cH_\Omega$ be an isometry of $(\Omega,d_\Omega)$ into a Hilbert space $\cH_\Omega$. Assume $\cH$ consists of $\cH_\Omega$-valued functions over $\cT$ and
\[  e(x)(t) = e_\Omega(x(t)),\hspace{4mm}\forall t\in\cT,x\in\cX. \]
Then $C$ given by~(\ref{cov2}) is a {\em kernel-valued kernel} over $\cT\times\cT$, given for $s,t\in\cT$ by
\begin{align*}
C(s,t) 
&= \frac1n\sum_{i=1}^n \big( e(X_i)(s)-\overline e(s)\big)\otimes \big(e(X_i)(t)-\overline e(t) \big) \nonumber \\
&= \frac1n\sum_{i=1}^n \big( e_\Omega(X_i(s))-\overline e(s)\big)\otimes \big(e_\Omega(X_i(t))-\overline e(t) \big).  \label{c2}
\end{align*}
Instead, DM's method is based on the metric covariance kernel $C_\text{DM}:\cT\times\cT\rightarrow\mR$,
\[  C_\text{DM}(s,t) = \frac1n\sum_{i=1}^n \big\langle e_\Omega(X_i(s))-\overline e(s), e_\Omega(X_i(t))-\overline e(t) \big\rangle . \]
We immediately have
\[  C_\text{DM}(s,t) = \tr(C(s,t)), \]
i.e., $C_\text{DM}$ entails some information loss relative to $C$.

The eigenvectors of $C$ are in the ``correct'' space of {\em $\cH_\Omega$-valued functions} over $\cT$ and kPCs are naturally defined by~(\ref{kpc}).
However, the eigenvectors of $C_\text{DM}$ are in the ``wrong'' space of {\em real-valued functions} over $\cT$, and it is not obvious how  metric-valued functional observations load on an eigenvector of $C_\text{DM}$. For this reason, DM needed to add some complexity with the new concepts named ``object FPCs'' and ``Fr\'echet scores'', which can be avoided with our approach. 






\bibliographystyle{elsarticle-harv}
\bibliography{stats}

\begin{thebibliography}{3}
\expandafter\ifx\csname natexlab\endcsname\relax\def\natexlab#1{#1}\fi
\providecommand{\url}[1]{\texttt{#1}}
\providecommand{\href}[2]{#2}
\providecommand{\path}[1]{#1}
\providecommand{\DOIprefix}{doi:}
\providecommand{\ArXivprefix}{arXiv:}
\providecommand{\URLprefix}{URL: }
\providecommand{\Pubmedprefix}{pmid:}
\providecommand{\doi}[1]{\href{http://dx.doi.org/#1}{\path{#1}}}
\providecommand{\Pubmed}[1]{\href{pmid:#1}{\path{#1}}}
\providecommand{\bibinfo}[2]{#2}
\ifx\xfnm\relax \def\xfnm[#1]{\unskip,\space#1}\fi
\bibitem[{Dubey and Mueller(2019)}]{dm19}
\bibinfo{author}{Dubey, P.}, \bibinfo{author}{Mueller, H.G.},
  \bibinfo{year}{2019}.
\newblock \bibinfo{title}{Functional models for time-varying random objects}.
\newblock \bibinfo{journal}{arXiv preprint arXiv:1907.10829} .
\bibitem[{{Human Mortality Database}(2019)}]{mortality19}
\bibinfo{author}{{Human Mortality Database}}, \bibinfo{year}{2019}.
\newblock \bibinfo{title}{{University of California, Berkeley (USA)} and {Max
  Planck Institute for Demographic Research (Germany)}}.
\newblock \bibinfo{note}{Available at \url{www.mortality.org} or
  \url{www.human-mortality.org}}.
\bibitem[{Sch{\"o}lkopf et~al.(1998)Sch{\"o}lkopf, Smola and
  M{\"u}ller}]{ssm98}
\bibinfo{author}{Sch{\"o}lkopf, B.}, \bibinfo{author}{Smola, A.},
  \bibinfo{author}{M{\"u}ller, K.R.}, \bibinfo{year}{1998}.
\newblock \bibinfo{title}{Nonlinear component analysis as a kernel eigenvalue
  problem}.
\newblock \bibinfo{journal}{Neural computation} \bibinfo{volume}{10},
  \bibinfo{pages}{1299--1319}.

\end{thebibliography}

\end{document}